\documentclass[12pt]{amsart}

\usepackage{amsmath}
\usepackage{amssymb}
\usepackage{latexsym}
\newtheorem{theorem}{Theorem}[section]
\newtheorem*{theorem*}{Theorem} 
\newtheorem*{lemmaA}{Lemma A}
\newtheorem{corollary}[theorem]{Corollary} 
\newtheorem{lemma}[theorem]{Lemma} 
\newtheorem{proposition}[theorem]{Proposition} 
\newtheorem{definition}[theorem]{Definition} 
 
\newenvironment{proof*}{\vskip 2mm\noindent {}}{\hfill$\square$ \vskip 2mm}
\numberwithin{equation}{section}

\newcommand{\C}{{\mathbb{C}}}
\newcommand{\Dd}{{\mathbb{D}}}
\newcommand{\Dbb}{{\mathbb{D}}}

\newcommand{\Zz}{{\mathbb{Z}}}

\renewcommand{\a}{\alpha}
 
\renewcommand{\d}{\delta} 
 
\newcommand{\e}{\varepsilon}
\newcommand{\eps}{\varepsilon}
\newcommand{\f}{\varphi} 
\newcommand{\g}{\gamma} 

\renewcommand{\L}{\Lambda}

\renewcommand{\t}{\theta} 

\newcommand{\calf}{{\mathcal F}}

\newcommand{\eit}{e^{i\theta}}
\newcommand{\dm}{\overline d_m}
\newcommand{\diz}{(1-|z|)}
\newcommand{\diszk}{(1-|z_k|)}
\newcommand{\szk}{\{z_k\}}


\thanks{Most of this paper was written while  the first-named author
was visiting  Universit\'e Bordeaux I in the Fall of 2002, supported by CNRS.
The first-named author thanks these institutions for their hospitality and
 support, and is especially grateful to
Alexander Borichev who arranged this visit.
Research of the first-named author is partially supported also by the Russian
Foundation of Basic Research, Grant 01-01-00608, and by the Ministry of
Education of Russia, Grant 00-1.0-199.}

\begin{document} 

\title[Summatory conditions for bounded functions]{Equivalence of summatory conditions along sequences
for bounded holomorphic functions}

\author{Vladimir Ya. Eiderman \ \ \ Pascal J. Thomas} 

\begin{abstract}
A sequence of points $z_k$ in the unit disk is said to be
thin for a given positive function $\rho$, if there is a nontrivial bounded
holomorphic function such that the infinite series 
$\sum_k \rho(1-|z_k|)|f(z_k)|$ converges. All sequences will be assumed
hyperbolically separated. We give necessary and sufficient 
conditions for the problem of thinness of a sequence to be non-trivial (one
 way or the other), and for two different positive functions $\rho_1, 
 \rho_2$ to 
give rise to the same thin sequences. Along the way, some concrete
conditions (necessary or sufficient) for a sequence to be thin are obtained. 
\end{abstract}

\maketitle

{\it This paper is dedicated to the memory of Matts Ess\'en.}
\vskip1cm

\section{Introduction and main results} 
\label{intro}

\footnotetext{{\it 2000 Mathematics Subject Classification} 30D50}

Let $\calf$ be a certain class of holomorphic functions in a domain
$G\subset\C$. N.~Nikolski \cite{Ni} and S.~Khavinson \cite{Kh} established the
existence of positive numbers $\{t_k\}$ (satisfying certain additional
conditions) and of a sequence of points $\{z_k\}\subset G$ such that
if $\sum_k t_k|f(z_k)|<\infty$, then $f\equiv0$ for every $f\in\mathcal F$.
Theorems of such kind arise in various applications: operator theory
\cite{Ni}, approximation theory \cite{Kh}, \cite{E3}, sampling, etc. They are also of
independent interest.

We shall consider the case when $G$ is the unit
disk $\Dbb=\{|z|<1\}$, $|z_k|\to1$ as $k\to\infty$,
and $\mathcal F$ is the class $H^\infty(\Dbb)$ of
bounded holomorphic functions in $\Dbb$. Note that since functions in the
Nevanlinna class (also called ``of bounded characteristic") can be written
as the quotient of two bounded functions, all our results also apply
automatically to that class. In \cite{CTW}, \cite{E1}, \cite{E2}, \cite{EE}
 sufficient
conditions on $\{z_k\}$ were obtained in the case $t_k\equiv1$. In
\cite{NPT}, \cite{EE} the problem of description of sequences 
$\{z_k\}$ was studied
for $t_k=1-|z_k|$ (see \cite{EE} for a more detailed survey).

In order to avoid summing repeatedly values of the function $f$ at nearby
points, we impose a uniform condition of discreteness on our sequence,
namely that it be separated in the pseudohyperbolic distance $d_G$.

\begin{definition} 
\label{separ}
We shall say that $\{z_k\}$ is \emph{separated}, if
\begin{equation}
\label{eqsepar}
\inf_{n\ne k} d_G (z_j, z_k) :=
\inf_{n\ne k}\biggl|\frac{z_n-z_k}{1-\overline{z}_nz_k}\biggr|=\d>0.
\end{equation}
\end{definition} 

A natural special case of
the problem is that of radial coefficients:
   $t_k=\rho(1-|z_k|)$, where $\rho(t),\ t\in
(0,\,1)$, is a positive function. It will be convenient to write the function
$\rho(t)$ in the form $\rho(t)=\rho_\t(t)=te^{\t(t)}$.

\begin{definition} 
\label{defthin}
 A sequence $\{z_k\}$ is said to be {\it $\t$-thin}, if
there exists a function $f\in H^\infty$ such that $f\not\equiv0$ and
\begin{equation}
\label{eqthin}
\sum_k \rho_\t (1-|z_k|)|f(z_k)|
=\sum_k (1-|z_k|)e^{\t(1-|z_k|) }|f(z_k)|<\infty.
\end{equation}
A non-$\t$-thin sequence is said to be {\it $\t$-thick}. The class of $\t$-thin
separated sequences is denoted by $\L_\t$.
\end{definition} 

Obviously, $\L_{\t_1}\subset\L_{\t_2}$ for $\t_1>\t_2$.

Our main results (Theorem \ref{equivthin} and Proposition 
\ref{equivthinsmall}) determine which
pairs of functions $\t_1$, $\t_2$ give rise to the same class of
$\t$-thin separated sequences.
Those results exhibit a certain
hierarchy among separated sequences in the disk.
In particular, the classes $\L_\t$ corresponding to
$\rho(t) = t^\a$, $0\le \a <1$, are all the same, but the case
$\rho(t) = t$ is distinct. Before the beginning of the spectrum
(for sufficiently big $\rho$, in particular all such $\rho$ must verify
$\limsup_{t\to0}\rho(t)=\infty$) all separated non Blaschke sequences are
thick (Theorem \ref{existthin}).
After the other end of the spectrum (for "small" $\rho$, in particular
$\liminf_{t\to0} \rho(t)/t =0$ for all such $\rho$)
all separated sequences (even maximal separated sequences, or "nets")
are thin (Proposition \ref{existthick}).

Our first statement is about the set of all possible positive functions $\t$.

\begin{theorem} 
\label{existthin}
Suppose that $\t(t)$ is a positive continuous
nonincreasing function. There exists a non Blaschke
separated $\t$-thin sequence $\{z_k\}$ if and only if
\begin{equation}
\label{eqexist}
\int_0\frac{dt}{t\t(t)}=\infty.
\end{equation}
\end{theorem} 
Obviously, Theorem \ref{existthin}
 implies the existence of $\t$-thin sequences for
all negative functions $\t(t)$.

\begin{theorem} 
\label{equivthin}
 a) If $t_1, t_2$ are positive functions and
$\t_1(t)\approx\t_2(t)$, then
$\L_{\t_1}=\L_{\t_2}$.

b) If $\t_1(t)/\t_2(t)\to\infty$ as $t\to0$ and $\t_2(t)$ is a positive
continuous  nonincreasing function satisfying \eqref{eqexist},
then $\L_{\t_1}\ne\L_{\t_2}$, i.~e. there exists a sequence
$\{z_k\}$ such that $\{z_k\}\in\L_{\t_2}$, but $\{z_k\}\not\in\L_{\t_1}$.
\end{theorem} 

That \eqref{eqexist} is a sufficient condition for the existence of elements in
   $\L_\t$ will follow from part (b) of this result, applied to $\t_2=\t$.
The necessity part in Theorem \ref{existthin}
 will be proved in Section \ref{suffnec}, Theorem \ref{equivthin}
will be proved in Section \ref{proofs}.

Consider now functions $\t$ bounded above (including negative ones).

\begin{proposition} 
\label{existthick} 
Assume that $\rho_\t(t)$ is a continuous
nondecreasing function such that $\rho_\t(t)\le Ct$ with $C>0$. There
exists a separated $\t$-thick sequence $\{z_k\}$ if and only if
\begin{equation}
\label{eqexthick}
\int_0\frac{\rho_\t(t)\,dt}{t^2}=\infty,\quad \rho_\t(t)=te^{\t(t)}.
\end{equation}
\end{proposition} 

Proposition \ref{existthick}  implies the existence of $\t$-thick sequences for
all nonnegative functions $\t(t)$.

For simplicity of notation, we write $\rho_1 = \rho_{\t_1}$,
$\rho_2 = \rho_{\t_2}$.

\begin{proposition} 
\label{equivthinsmall} 
Assume that $\rho_1(t)$ satisfies the conditions
of Proposition \ref{existthick}, \eqref{eqexthick} holds and
$\rho_1(t)/\rho_2(t)\to\infty$ as $t\to0$. Then $\L_{\t_1}\ne\L_{\t_2}$.
Moreover, there exists a sequence
$\{z_k\}$ such that $\{z_k\}\not\in\L_{\t_1}$, but (1.2) with
$\rho_\t=\rho_2$ holds for every $f\in H^\infty$.
\end{proposition} 

Notice that the sufficient part of Proposition \ref{existthick} 
(the existence of a
$\t$-thick sequence when \eqref{eqexthick}
is satisfied) follows from Proposition \ref{equivthinsmall},
taking $\t_1=\t$ and an arbitrary $\t_2$.

In Section \ref{suffnec} we establish some preliminary results. 
Some of them are of
independent interest. For example, Proposition \ref{slices} gives conditions on
$\{z_k\}$ which imply that $\{z_k\}$ is $\t$-thick, as well as a sufficient
condition of existence of $\t$-thin sequences. We develop these results
in Section \ref{furtherprop}.

\section{Some sufficient and some necessary conditions}
\label{suffnec}

First, we would like to show that our problem,
which is a priori given in terms of infinite sums, reduces in many cases to a
problem about the rate of decrease of $|f(z)|$ as a function of $\diz$,
when $z$ is restricted to the sequence $\szk$.  Such problems have already
been investigated in \cite{PT} and numerous other works, see the
references in \cite{EE}.

\begin{proposition} 
\label{equivdec}
(a) Suppose that $\szk \in \L_\t$, $\theta(t) \ge 0$.
Then there exists a non identically vanishing
function $f_1\in H^\infty(\Dd)$ such that
$$
    |f_1(z_k)| \le  e^{-\theta \diszk} ,
$$
for all $k \in \Zz_+$.

(b) Conversely, suppose that
$$
\liminf_{t\to0} \frac{\t(t)}{\log \log \frac1t}=: L >0,
$$
that  $\szk$ is separated, non-Blaschke, and
that there exists a function $f\in H^\infty(\Dd)$ such that
$$
    |f(z_k)| \le  e^{-\theta \diszk} ,
$$
for all $k \in \Zz_+$. Then $\szk  \in \L_\t$.
\end{proposition} 

\begin{proof}
(a) There exists a non identically vanishing
function $f\in H^\infty(\Dd)$ such that
$\sum_k \diszk e^{\theta \diszk} |f(z_k)| < \infty$.
Define
$$
\{z'_j\} := \{ z_k :  |f(z_k)| \le e^{-\theta \diszk}  \}.
$$
Then
$$
\sum_{k:z_k \notin \{z'_j\}} \diszk
\le
\sum_{k:z_k \notin \{z'_j\}} \diszk e^{\theta \diszk}  |f(z_k)| <
\infty.
$$
So $\{ z_k \} \setminus \{z'_j\}$ is a Blaschke sequence.
Let $B$ denote the Blaschke
product with zeroes on $\{ z_k \} \setminus \{z'_j\}$,
then $f_1 := f B$ clearly
satisfies the desired conclusion. 

(b) Let
\begin{equation}
\label{eqslices}
Y_n:=\{z:1-2^{-n}\le|z|<1-2^{-n-1}\}.
\end{equation}
Since $\szk$ is separated, $\sum_{k:z_k \in Y_n} \diszk \le C_\delta$. Let $m
\in \Zz_+$ to be chosen later, $f_1 := f^m$. Then, because of the above and
of the hypothesis, for any $\eps>0$,
\begin{multline*}
\sum_k \diszk e^{\theta \diszk} |f(z_k)^m|
\le
\sum_n \sum_{k:z_k \in Y_n} \diszk e^{-(m-1) \theta \diszk}
\\
\le
C_\eps + C_\delta \sum_n  e^{-(m-1) (L-\eps) \log n}
=
C_\eps + C_\delta \sum_n \frac{1}{n^{(m-1) (L-\eps)}} < \infty
\end{multline*}
for $m$ large enough (for $L=\infty$, replace $L-\e$ by any positive
number).
\end{proof}

\begin{proof*}{\it Proof of the necessity part in Theorem \ref{existthin}.}

The existence of a non-trivial $f_1\in H^\infty(\Dd)$ such that
$|f_1(z)| \le  e^{-\theta \diz}$ over the sequence $\szk$, which is
non-Blaschke and separated, implies that $g(t):=  e^{-\theta (1-t)}$ is
{\it not\/}  an
{\it essential minorant\/} in the language of \cite{LS}. The main result
of that paper is that $g$ is an essential minorant if and only if
$$
\int^1 \frac{dt}{(1-t)\log(\frac1{g(t)})} < \infty,
$$
which is exactly the negation of \eqref{eqexist}. 
Therefore, if \eqref{eqexist} does not
hold, there are no sequences in $\L_\t$. 
\end{proof*}

For a given sequence $\{z_i\}$ let $N_m$ be a number of points $z_i$ in
$Y_m$, where $Y_m$ is as in \eqref{eqslices}.

\begin{lemma}
\label{excepindices}
Let $\sum_{m}N_m2^{-m}=\infty$. Then for every $f\in H^\infty$,
$f\not\equiv0$, there exists a set $\mathcal J$ of indices with the
following properties: $\sum_{m\in\mathcal J}N_m2^{-m}<\infty$, and
$$
|f(z_k)|>\exp\{-C(f,\d)2^m/N_m\} \quad\text{for at least $N_m/2$
points } z_k\in Y_m,\ m\not\in\mathcal J,
$$
where $C(f,\d)$ is a positive constant and $\d$ is the constant from 
\eqref{eqsepar}.
\end{lemma}

Lemma \ref{excepindices} was in fact proved in \cite{EE} 
(see proof of Proposition 7.9 (a)),
although it was not formulated in this explicit form. In Section 4 we give
another, essentially self-contained proof (without
resorting to Govorov and Grishin as in \cite{EE}).

\begin{proposition} 
\label{slices}
Let $\t(t)$ be a positive nonincreasing
function and $\rho_\t(t)=te^{\t(t)}$.

a) Suppose that $\{z_k\}$ is a separated sequence and
\begin{equation}
\label{sumdiv}
\sum_{m\not\in\mathcal J}N_m\rho_\t(2^{-m})\exp\{-\g2^m/N_m\}=\infty
\quad\text{for every }\g>0
\end{equation}
and for every set $\mathcal J$ of indices with
$\sum_{m\in\mathcal J}N_m2^{-m}<\infty$. Then $\{z_k\}\not\in\L_\t$.

b) Suppose that a set $\mathcal L$ of indices and a sequence $\{N_m\}$ of
nonnegative integers are such that $N_m\uparrow\infty$,
$N_m2^{-m}\downarrow0$ as $m\in\mathcal L$, $m\to\infty$, and
\begin{equation}
\label{sumconv}
\sum_{m\in\mathcal L} N_m\rho_\t(2^{-m})\exp\{-\g2^m/N_m\}<\infty
\quad\text{for some }\g>0.
\end{equation}
Then there exists a separated $\t$-thin sequence $\{z_k\}$, for which $N_m$
are the given numbers.
\end{proposition} 

\vskip.4cm

{\bf Examples} 

a) If $N_m=p_m2^m/(\t(2^{-m})+\log m)$,
where $p_m\to\infty$, then $\{z_k\}\not\in\L_\t$.

b) If $N_m=p2^m/(\t(2^{-m})+\log m)$ with a positive
constant $p$, then there exists a separated sequence $\{z_k\}\in\L_\t$,
for which $N_m$ are the given numbers.

Proposition~\ref{slices} generalizes and refines Proposition 7.9 in \cite{EE}.
Its proof will be given in Section \ref{prooflemma}.

\section{Proofs of the main results}
\label{proofs}

\begin{proof*}{\it Proof of Theorem \ref{equivthin} (a).}
It will be enough to show that if $\theta_1 \ge c \theta_2$, for some
$c>0$, then $\Lambda_{\t_1} \subset \Lambda_{\t_2}$. Without loss of
generality, assume that $c<1$.

Suppose that $\szk \in \Lambda_{\t_1}$. By Proposition \ref{equivdec}, 
we can pick a non
identically vanishing function $f_1 \in H^\infty(\Dd)$ such that
$$
|f_1(z_k)| \le e^{-\theta_1 \diszk} , \quad \forall k \in \Zz_+.
$$
Pick an integer $m \ge 1+ c^{-1}$. Then
$$
|f_1 (z_k)^m| \le
    e^{-c^{-1}\theta_1 \diszk}  e^{\theta_1 \diszk}|f_1 (z_k)|,
$$
so
\begin{eqnarray*}
\sum_k e^{\theta_2 \diszk} |f_1 (z_k)^m|
&\le&
\sum_k e^{\theta_2 \diszk-c^{-1}\theta_1 \diszk}
e^{\theta_1 \diszk}|f_1 (z_k)|
\\
&\le&
\sum_k  e^{\theta_1 \diszk}|f_1 (z_k)| < \infty,
\end{eqnarray*}
so $f_1^m$ is a non
identically vanishing function in $H^\infty(\Dd)$ which proves that
$\szk \in \Lambda_{\t_2}$. 

\end{proof*}

\begin{proof*}{\it Proof of Theorem \ref{equivthin} (b).} Assume that 
\eqref{eqexist} holds with $\t=\t_2$.
First we construct a sequence $\szk \in \L_{\t_2}$.

We note that \eqref{eqexist} is equivalent to the following condition:
\begin{equation}
\label{eqexistser}
\sum_{m=1}^\infty[\t_2(2^{-m})]^{-1}=\infty.
\end{equation}
Without loss of generality we may assume that $\t_2(t)\uparrow\infty$ as
$t\to0$. We show that there exist a set $\mathcal L$ of indices and a
nondecreasing sequence $\{N_m\}$ of positive integers
with the following properties:
\begin{equation}
\label{eqn32}
\sum_{m\in\mathcal L}[\t_2(2^{-m})]^{-1}=\infty,\quad
\sum_{m\in\mathcal L}[\t_2(2^{-m})]^{-2}<\infty,
\end{equation}
\begin{equation}
\label{eq33}
2^{m-1}/\t_2(2^{-m})<N_m\le 2^{m}/\t_2(2^{-m}),\quad
N_m2^{-m}\downarrow0 \text{ as } m\in\mathcal L,\ m\to\infty.
\end{equation}
Denote $\e_m=[\t_2(2^{-m})]^{-1}$. Then $\e_m\downarrow0$ and
$\sum_m\e_m=\infty$. Choose integers $m_n$, $M_n$ such that $m_n<M_n\le 
m_{n+1}$ and
$$
\frac1{2n}<\sum_{m=m_n+1}^{M_n}\e_m<\frac1{n},
$$
and set $\mathcal L_1=\bigcup_n\{m_n+1,\dots,M_{n}\}$. Then \eqref{eqn32} holds with
$\mathcal L=\mathcal L_1$. Let
$$
\mathcal J_m=\{j:j>m,\ 2^j\e_j\le 2^m\e_m\},\quad
\mathcal L=\mathcal L_1\setminus\bigcup_{m\in\mathcal L_1}\mathcal J_m.
$$
Since
$$
\sum_{j\in\mathcal J_m}\e_j\le2^m\e_m\sum_{j=m+1}^\infty2^{-j}=\e_m,
$$
the relations \eqref{eqn32} with $\mathcal L=\mathcal L_1$ 
imply the same relations with
$\mathcal L$. For $m,k\in\mathcal L$, $m<k$ we have $2^{m-k}\e_m<\e_k<\e_m$. Fix
$m_1\in\mathcal L$ and set $N_{m_1}=[2^{m_1}\e_{m_1}]$, where $[\cdot]$
denotes the integral part of a number. Suppose that $N_m$ is already defined,
and $k=\min\{j:j\in\mathcal L,\ j>m\}$. Then we set
$$
N_k=\min\{2^{k-m}N_m,\,[2^k\e_k]\}.
$$
Since $2^m\e_m<2^k\e_k$, we have
$N_k\ge N_m$.

Now we prove \eqref{eq33}. If $N_k=[2^k\e_k]$, then the inequalities
in \eqref{eq33} are evident. If $N_k=2^{k-m}N_m$, then by induction we have
$$
2^{k-1}\e_k<2^{k-1}\e_m=2^{k-m}2^{m-1}\e_m<2^{k-m}N_m=N_k.
$$
Moreover, the inequality $N_k\le N_m2^{k-m}$ implies that
$N_k2^{-k}\le N_m2^{-m}$. Thus, \eqref{eq33} is established.

Let $l_m=N_m2^{-m}$.
By \eqref{eq33} we have $\t_2(2^{-m})\le1/l_m<2\t_2(2^{-m})$, and
\begin{eqnarray*}
N_m\rho_2(2^{-m})\exp\{-\g2^m/N_m\}&=&
l_m\exp\{\t_2(2^{-m})-\g/l_m\}\\
&\le& l_m\exp\{-1/l_m\}<l_m^2,
\end{eqnarray*}
if $\g$ is sufficiently big. By \eqref{eqn32}, $\sum_{m\in\mathcal L}l_m^2<\infty$.
Thus, the condition \eqref{sumconv} is satisfied, and the existence of
a separated $\t_2$-thin sequence $\{z_k\}$ follows from 
Proposition~\ref{slices}~b).

Now we show that $\{z_k\}\not\in\L_{\t_1}$. Let
$p_m=\t_1(2^{-m})/\t_2(2^{-m})$. By \eqref{eq33},
$\t_2(2^{-m})\le 2^m/N_m<2\t_2(2^{-m})$. For $\rho_1(t)=te^{\t_1(t)}$
and for a fixed constant $\g$ we have
\begin{eqnarray*}
N_m\rho_1(2^{-m})\exp\{-\g2^m/N_m\}&=&
l_m\exp\{\t_1(2^{-m})-\g2^m/N_m\}\\
&\ge& l_m\exp\{\t_2(2^{-m})(p_m-2\g)\}>l_m,
\end{eqnarray*}
for sufficiently big $m\in\mathcal L$. By \eqref{eqn32} and \eqref{eq33}
$\sum_{m\in\mathcal L}l_m=\infty$. We see that \eqref{sumdiv} holds, and by
Proposition~\ref{slices} a) $\{z_k\}\not\in\L_{\t_1}$. 
\end{proof*}

\begin{proof*}{\it Proof of the necessity part of Proposition \ref{existthick}.}

Suppose that \eqref{eqexthick} fails. By taking a Riemann sum, we will show that
any separated sequence is thin, even
better : that we can take as function $f$ the constant $1$.
The separation condition \eqref{eqsepar} implies that every
Whitney square contains at most $M=M(\d)$ points $z_k$. Hence, $N_m=O(2^m)$
for every separated sequence $\{z_k\}$.
Then
$$
\sum_k\rho\diszk\le\sum_m N_m\rho(2^{-m})\le
C\sum_m2^m\rho(2^{-m});
$$
the convergence of the last series is equivalent
to the negation of \eqref{eqexthick}.
\end{proof*}

To prove Proposition \ref{equivthinsmall},
we need the following elementary fact about numerical series.

\begin{lemma} 
\label{series}
If $\sum_n 2^n\rho_1(2^{-n})=\infty$,
$\rho_1(t)\le Ct$ and $\rho_1(t)/\rho_2(t)\to\infty$ as $t\to0$,
we can choose a subset of the integers $\mathcal A \subset \Zz_+$ such that
\begin{equation}
\label{rhoconv}
\sum_{n\in\mathcal A} 2^n \rho_2 (2^{-n}) < \infty,
\end{equation}
and
\begin{equation}
\label{rhodiv}
\sum_{n\in\mathcal A} 2^n \rho_1 (2^{-n}) = \infty.
\end{equation}
\end{lemma}

\begin{proof}
We can recursively define a strictly
increasing sequence of integers $\{n_j\}$ such that
\begin{enumerate}
\item $\rho_2(2^{-n}) \le 2^{-j} \rho_1(2^{-n})$, $n \ge n_j$;
\item $\sum_{n_j \le n <n_{j+1}} 2^n \rho_1(2^{-n}) \ge 1$.
\end{enumerate}
Now we pick
$$
n'_j := \min \{m \ge n_j :
\sum_{n_j \le n \le m} 2^n \rho_1(2^{-n}) \ge 1 \},
$$
and set
$$
\mathcal A_j:= \{ n \in \Zz_+ : n_j \le n \le n'_j \},
\quad \mathcal A := \cup_j\mathcal A_j.
$$
It is clear
that $\sum_{\mathcal A} 2^n \rho_1(2^{-n}) \ge \sum_j 1 = \infty$.
Because $2^n \rho_1(2^{-n}) \le C$ for all $n$,
$\sum_{\mathcal A_j} 2^n \rho_1(2^{-n}) \le 1 +C$ for all $j$, therefore by (1)
$$
\sum_{\mathcal A} 2^n \rho_2(2^{-n}) \le \sum_j\sum_{\mathcal A_j} 2^n\rho_2(2^{-n})
\le \sum_j 2^{-j} (1+C) < \infty.
$$
\end{proof}

\begin{proof*}{\it Proof of Proposition \ref{equivthinsmall}.} 
Consider
the sequence $S=\{z_{m,j}: m \in\mathcal A,\ 0 \le j < 2^m\}$, where
$$
z_{m,j} := (1-2^{-m})\exp(2^{-m}2\pi ij), \quad 0 \le j < 2^m,
$$
and $\mathcal A$ is the set from Lemma \ref{series}. For any $f\in H^\infty$, by 
\eqref{rhoconv},
$$
\sum_{z\in S} \rho_2\diz |f(z)| \le \|f\|_\infty \sum_{z\in S} \rho_2\diz
=  \|f\|_\infty \sum_{n\in\mathcal A} 2^n \rho_2 (2^{-n}) < \infty.
$$

Now we will show that $S \notin \L_{\t_1}$.
Let $f\in H^\infty$, $f\not\equiv0$, be given, and let $\mathcal J$ be the set
from Lemma~\ref{excepindices}. Since $N_m=2^m$ as $m\in\mathcal A$, the
set $\mathcal J$ is finite, and by \eqref{rhodiv} we have
$$
\sum_k \rho_1(1-|z_k|)|f(z_k)|
\ge\frac12\sum_{m\in\mathcal A\setminus\mathcal J}
2^m\rho_1(2^{-m})\exp\{-C(f,\d)\}=\infty.
$$
\end{proof*}

\section{Proof of Lemma \ref{excepindices} and Proposition~\ref{slices}.}
\label{prooflemma}

Our alternative prove of Lemma \ref{excepindices} is based on the following
Lemma, which is proved in \cite[p.~124, lines 3 to 17]{NPT}
and is the main ingredient in \cite[Lemma~3]{PT}.

\begin{lemmaA}
Let $f\in H^\infty(\Dd)$ and $\d' \in (0,1)$.
Then there exists a function $h$, positive and
harmonic on $\Dd$, such that for all $z$ for which $d_G(z,f^{-1}(0)) \ge
\delta'$,
$$
e^{-h(z)} \le |f(z)|.
$$
In particular, if $\szk$ is a separated sequence,
there exists a Blaschke sequence $b \subset \szk$ such that for any
$z_k \notin b$,
$$
e^{-h(z_k)} \le |f(z_k)|.
$$
\end{lemmaA}

\begin{proof*}{\it Proof of Lemma \ref{excepindices}.} By Lemma A, there exists 
a positive
harmonic function $h$ and a Blaschke sequence $b \subset \{z_{k}\}$
such that
$$
\mbox{For } z_k \notin b,\,  e^{-h(z_k)} < |f(z_k)|.
$$
Let
$$
\mathcal J:= \{m : \#b\cap Y_m \ge N_m/4 \}.
$$
Since $b$ is a Blaschke
sequence, we have $\sum_{m\in\mathcal J }N_m2^{-m}<\infty$.

Let $C(f,\d)>0$ be a constant to be chosen later. For a point $a\in Y_m$,
let $I_a \subset \partial \Dbb$ be the arc centered at
$a/|a|$, of length $2^{-m}$.
By the separation
condition, the $I_{z_k}$ for $z_k \in Y_m$ form a covering of finite
multiplicity of $\partial \Dbb$, that is to say,
there exists a positive constant $C_1(\d)$ such that
$$
\sum_{z_k \in Y_m} \chi_{I_{z_k}} \le C_1(\d).
$$
By Harnack's inequality, for $\theta$ such that
$\eit z_k/|z_k| \in I_{z_k}$, we have
$$
h(\eit (1-2^{-m}) z_k/|z_k|) \ge C_2 h(z_k).
$$
For $m\in \mathcal J$, let
$$
N'_m:= \#\left\{ z_k \in Y_m \setminus b : h(z_k) > C(f,\d) 2^m/N_m
\right\}.
$$
The mean-value property for harmonic functions gives
$$
h(0) = \int_0^{2\pi} h( (1-2^{-m})\eit ) \frac{d\theta}{2\pi}
\ge C_1(\d)^{-1} \cdot N'_m \cdot 2^{-m} \cdot C_2 C(f,\d) 2^m/N_m,
$$
therefore
$$
N'_m \le (C_2 C(f,\d))^{-1} C_1(\d) h(0) N_m < N_m/4
$$
for
$C(f,\d)$ large enough. At the remaining points $z_k$, of which there are at
least $N_m/2$, we do have
$$
|f(z_k)| \ge e^{-h(z_k)}  \ge \exp\{-C(f,\d)2^m/N_m\}. \qquad 
$$
\end{proof*}

\begin{proof*}{\it Proof of Proposition \ref{slices} (a).} 
Under conditions of Proposition \ref{slices},
we have
$$
\sum_{m}N_m2^{-m}=\infty.
$$
Moreover, $\rho_\t(1-|z|)=(1-|z|)e^{\t(1-|z|)}>
2^{-m-1}e^{\t(2^{-m})}=\frac12\rho_\t(2^{-m})$, $z\in Y_m$. Let
$f\in H^\infty$, $f\not\equiv0$, be given, and let $\mathcal J$ be the set
from Lemma~\ref{excepindices}. If $\g>C(f,\d)$, then by \eqref{sumdiv}
\begin{eqnarray*}
\sum_k \rho_\t(1-|z_k|)|f(z_k)|&
\ge &\frac12\sum_{m\not\in\mathcal J}\sum_{z_k\in Y_m}\rho_\t(2^{-m})|f(z_k)|\\
&\ge& \frac14\sum_{m\not\in\mathcal J}
N_m\rho_\t(2^{-m})\exp\{-\g2^m/N_m\}=\infty.
\end{eqnarray*}
\end{proof*}

\begin{lemma} 
\label{lemmaEE}
Suppose that $P_1(t),\ t\in(0,1)$, is a
nondecreasing function, $P_1(t)(1-t)$ is a nonincreasing function
tending to $0$ as $t\to1$, and $d_m \le 1$ are given positive numbers. Then
there are a sequence $\{z_k\}\subset\Dbb$ and an absolute constant $M$
with the following properties: all points $z_k$ are situated on the
circumferences $|z|=1-2^{-m}$, $m=1,2,\dots$;
$$
\min\{|z_j-z_k|: j\ne k,\ |z_j|=|z_k|=1-2^{-m}\}\ge d_m;
$$
\begin{equation}
\label{eq41}
N_z(1-2^{-m})>M(d_mP_1(1-2^{-m}))^{-1}+1,\quad m=1,2,\ldots,
\end{equation}
where $N_z(1-2^{-m})$ is the number of points $z_k$ with $|z_k|=1-2^{-m}$,
and for every $T>0$ there is a nontrivial function $f\in H^\infty$ such
that
\begin{equation}
\label{logfdec}
\log|f(z_k)|\le-P_1(|z_k|)T,\quad k=1,\,2\,\ldots
\end{equation}
\end{lemma}

We omit the proof, which is exactly the same as the proof of Lemma~4.3,
part 2 in \cite{EE} with $\phi(t_j)=1/d_j$ and $t_j=1-2^{-j}$.
It follows from Lemma 4.3 in \cite{EE} that Lemma \ref{lemmaEE} is sharp.

\begin{proof*}{\it Proof of Proposition \ref{slices} (b).} Pick $d_m=2^{-m}$.
By assumption, $N_m2^{-m}\downarrow0$,
$\{N_m\}$ is a nondecreasing sequence and $N_m\to\infty$ on $\mathcal L$.
Hence, there exists a nondecreasing function $P_1(t)$ such that
$$
P_1(1-2^{-m})=M2^m/N_m,\quad m\in\mathcal L,
$$
where $M$ is the constant from \eqref{eq41}, $P_1(t)\uparrow\infty$,
$P_1(t)(1-t)\downarrow0$ as $t\to1$. Let $\{z_k\}$
be the sequence from Lemma~4.1. By \eqref{eq41},
$$
N_z(1-2^{-m})>N_m,\quad m\in\mathcal L.
$$
If we choose a sufficiently big number $T$, then by \eqref{logfdec}
$$
\log|f(z_k)|\le-M2^mT/N_m<-\g2^m/N_m,\quad |z_k|=1-2^{-m},\quad m\in\mathcal L.
$$
For every $m\in\mathcal L$ we choose $N_m$ points from $\{z_k\}$ with
$|z_k|=1-2^{-m}$.  We also denote this subsequence by $\{z_k\}$.
Using \eqref{sumconv}, we have
\begin{eqnarray*}
\sum_k\rho_\t(1-|z_k|)|f(z_k)|
&=& \sum_{m\in\mathcal L}\sum_{z_k\in Y_m}\rho_\t(1-|z_k|)|f(z_k)|\le\\
&\le& c\sum_{m\in\mathcal L}N_m\rho_\t(2^{-m})\exp\{-\g2^m/N_m\}<\infty.
\end{eqnarray*}
Thus, $\{z_k\}\in\L_\t$, and the proof of Proposition~2.3 is complete.
\end{proof*}

\section{Further properties of thin sequences}
\label{furtherprop}

The following assertion was proved in \cite{EE} (see Proposition~7.8).
{\it Let $N_m$ be given nonnegative integers, and let
$$
z_{m,k}=(1-2^{-m})\exp(2\pi ik/N_m),\quad 1\le k\le N_m,\quad m=1,2,\dots
$$
for $N_m>0$. Then $\{z_{m,k}\}\not\in\L_\t$ for $\rho_\t\equiv1$ if and
only if $\{z_{m,k}\}$ is a non-Blaschke sequence, that is
$\sum_mN_m2^{-m}=\infty$.}

On the other hand, if we know only, that $\{z_{m,k}\}$ is a separated
sequence, we need essentially stronger conditions on $\{N_m\}$ (see
Proposition~\ref{slices} with $\rho_\t\equiv1$). 
Thus, there is a certain connection
between the dispersion of points over annuluses $Y_m$ and the massivity
of a sequence, i.~e. numbers $N_m$ (we also can interpret this relation as
a connection between the dispersion of points and the possible decay of
nontrivial bounded function over a sequence). In the present section we
investigate this connection. In particular, we develop 
Lemma~\ref{excepindices},
Proposition~\ref{slices} and some results in \cite{EE}, and give another approach to
obtaining such assertions, which is based on estimates of subharmonic
functions outside exceptional sets.

Let a sequence $\{z_k\}$ be given. For every $m$, for which $Y_m$
contains at least six points $z_k$, and for every $z_k\in Y_m$ we set
$$
d_{m,k}:=\min\{|z_j-z_k| : j\ne k,\ z_j,z_k\in Y_m\}.
$$
Let $\dm$ be the mean value of $[N_m/6]$ smallest numbers $d_{m,k}$
(here $[\cdot]$ is the integral part of a number).
Our main tool is the following lemma.

\begin{lemma}
\label{slowdec}
Let $\{z_k\}$ be a separated sequence. Let a set $\mathcal K$ of indices
be such that $N_m\ge6$ as $m\in\mathcal K$, and
$\sum_{m\in\mathcal K}N_m2^{-m}=\infty$. Then for every $f\in H^\infty$,
$f\not\equiv0$, there exists a set $\mathcal J$ of indices with the
following properties: $\sum_{m\in\mathcal J}N_m2^{-m}<\infty$, and
\begin{equation}
\label{eqslowdec}
|f(z_k)|>\exp\{-C(f,\d)(N_m\dm)^{-1}\} 
\end{equation}
for at least $N_m/2$
points $z_k\in Y_m$, when $m\in\mathcal K\setminus\mathcal J$, 
where $C(f,\d)$ is a positive constant and $\d$ is the constant from 
\eqref{eqsepar}.
\end{lemma}

Since $\dm>\text{Const}\,2^{-m}$, Lemma \ref{slowdec} 
refines Lemma \ref{excepindices}. 
We shall prove Lemma~\ref{slowdec} in Section \ref{last}.

\begin{corollary}
\label{slicessep}
 Let $\t(t)$ be a positive nonincreasing
function and $\rho_\t(t)=te^{\t(t)}$.

a) Suppose that $\{z_k\}$ is a separated sequence. Let $\mathcal K$ be a set
of indices satisfying conditions of Lemma \ref{slowdec} and such that
\begin{equation}
\label{sumNdiv}
\sum_{m\in\mathcal K\setminus\mathcal J}N_m\rho_\t(2^{-m})
\exp\{-\g(N_m\dm)^{-1}\}=\infty \quad\text{for every }\g>0
\end{equation}
and for every set $\mathcal J$ of indices with
$\sum_{m\in\mathcal J}N_m2^{-m}<\infty$. Then $\{z_k\}\not\in\L_\t$.

b) Suppose that a set $\mathcal L$ of indices, a sequence $\{N_m\}$ of
nonnegative integers and a sequence $\{d_m\}$ of positive numbers
are such that $1\ge d_m\ge 2^{-m}$, $\{d_m2^m\}$ is a nondecreasing
sequence, $N_m\uparrow\infty$,
$N_md_m\downarrow0$ as $m\in\mathcal L$, $m\to\infty$, and
\begin{equation}
\label{sumNconv}
\sum_{m\in\mathcal L} N_m\rho_\t(2^{-m})\exp\{-\g(N_md_m)^{-1}\}<\infty
\quad\text{for some }\g>0.
\end{equation}
Then there exists a separated $\t$-thin sequence $\{z_k\}$ with
$|z_k|=1-2^{-m}$ as $z_k\in Y_m$, for which $N_m$ are the given numbers
and $\dm\ge d_m$.
\end{corollary}
\begin{proof} (a) The proof is essentially the same as the proof of
Proposition~\ref{slices} (a) with the obvious corrections: we use 
Lemma~\ref{slowdec}
instead of Lemma~\ref{excepindices}.

(b) We argue by analogy with the proof of Proposition~\ref{slices} (b).
Let $P_1(t)$ be a nondecreasing function such that $P_1(t)\uparrow\infty$,
$P_1(t)(1-t)\downarrow0$ as $t\to1$ and
$$
P_1(1-2^{-m})=M(N_md_m)^{-1},\quad m\in\mathcal L,
$$
where $M$ is the constant from \eqref{eq41}. Let $\{z_k\}$
be the sequence from Lemma~\ref{lemmaEE}. By \eqref{eq41},
$$
N_z(1-2^{-m})\ge N_m,\quad m\in\mathcal L.
$$
By \eqref{logfdec}, for sufficiently big numbers $T$ we have
$$
\log|f(z_k)|\le-MT(N_md_m)^{-1}<-\g(N_md_m)^{-1},
\quad |z_k|=1-2^{-m},\quad m\in\mathcal L.
$$
For every $m\in\mathcal L$ we choose $N_m$ points from $\{z_k\}$ with
$|z_k|=1-2^{-m}$.  The obtained sequence we also denote by $\{z_k\}$.
Using \eqref{sumNconv}, we have
\begin{eqnarray*}
\sum_k\rho_\t(1-|z_k|)|f(z_k)|
&=&\sum_{m\in\mathcal L}\sum_{z_k\in Y_m}\rho_\t(1-|z_k|)|f(z_k)|\le\\
&\le& c\sum_{m\in\mathcal L}N_m\rho_\t(2^{-m})\exp\{-\g(N_md_m)^{-1}\}
<\infty.
\end{eqnarray*}
\end{proof}

Corollary \ref{slicessep} shows the connection between the dispersion of points over
annuluses $Y_m$ and the massivity of a sequence $\{z_k\}$. It is clear
that Corollary~\ref{slicessep} is a direct generalization of 
Proposition~\ref{slices}.
Another special case gives the following assertion.

\begin{corollary} 
\label{slicesfix}
Let $N_m\dm\ge c>0$ as $m\in\mathcal K$ and
$\rho_\t(t)\ge c_1t,\ c_1>0$. Then a separated sequence  $\{z_k\}$ is
$\t$-thick if and only if
\begin{equation}
\label{eq54}
\sum_{m\in\mathcal K}N_m2^{-m}=\infty. 
\end{equation}
\end{corollary}
\begin{proof} Clearly, if $\{z_k\}\not\in\L_\t$, then \eqref{eq54} holds.

Conversely, the condition \eqref{eq54} and the inequality $\rho_\t(t)\ge c_1t$
imply \eqref{sumNdiv}, and our assertion follows from 
Corollary~\ref{slicessep}. 
\end{proof}

Corollary \ref{slicesfix} is a generalization of Proposition 7.8 in 
\cite{EE}, quoted
at the beginning of this section. It is easy to see 
that Corollary \ref{slicesfix}
is not correct for $\rho_\t(t)$ such that $\liminf_{t\to0}\rho_\t(t)/t=0$.
Indeed, in this case there exist numbers $t_j\in(0,1)$ for which
$t_{j+1}/t_j<1/2$ and $\sum_{j=1}^\infty\rho_\t(t_j)/t_j<\infty$.
Let $\mathcal K=\{m_j\}$, where $2^{-m_j-1}<t_j\le2^{-m_j}$, and let
$N_{m_j}=2^{m_j}$, $z_{m_j,k}=(1-t_j)\exp(2\pi ik/N_{m_j})$,
$1\le k<N_{m_j}$. Then $\{z_{m_j,k}\}$ is a separated $\t$-thin
sequence satisfying \eqref{eq54}.

\section{Proof of Lemma \ref{slowdec}.}
\label{last}

We need certain preliminary definitions and results.
Let $f$ be a meromorphic function in $\Dbb$, and let $\{a_k\}$ and
$\{b_k\}$ be zeros and poles of $f$, respectively. We assume that 0 is
neither a zero nor a pole of $f$. The Nevanlinna characteristic of $f$
is defined by (cf. for example Hayman \cite{Ha}, p.~4)
$$
T_f(r)=\frac1{2\pi}\int_0^{2\pi}\log^+|f(re^{i\f})|\,d\f+
\sum_{|b_i|<r}\log\frac r{|b_i|}.
$$
Since $T_f(r)$ is a nondecreasing function of $r$,
$\lim\limits_{r\uparrow1}T_f(r)=:T_f(1)$ exists for functions $f$ in
the Nevanlinna class $N$, i.~e., functions of bounded characteristic.

To $f\in H^\infty(\Dbb)$ with zeros  $\{a_n\},\ a_n\ne0$, we associate a
non-negative measure $\eta$ on the closed unit disc $\overline \Dbb$
defined as follows: for every set  $E\subset\overline\Dbb$ for which
$E\cap\partial\Dbb$ is measurable, we define
\begin{equation}
\label{caracE}
\eta(E)=\lim_{r\to1}\frac1{2\pi}\int_{e^{i\t}\in E}\log^-|f(re^{i\t})|\,d\t
+\sum_{a_n\in E}\log\frac1{|a_n|}.
\end{equation}
We note that $\eta(\overline\Dbb)=T_{1/f}(1)$. For the proof that this
limit exists, we refer to \cite[Appendix~1]{EE}.
We quote the following theorem in the form given
in \cite{EE} (see Theorem~2.10 in \cite{EE} and the remark after it).

\begin{theorem} (Govorov--Grishin) 
\label{GovGrish}
Suppose that $f\in
H^\infty(\Dbb)$, $f(z)\not\equiv0$. For every $P>1$ there exists a system
of discs $D^k(w_k,r_k)$ such that
\begin{equation}
\label{logfP}
\log|f(z)|>-A_1PT_{1/f}(1),\quad z\in\Dbb\setminus\bigcup_kD^k,
\end{equation}
\begin{equation}
\label{sumradii}
\sum_kr_k<1/P,
\end{equation}
and every $z\in\bigcup_kD^k$ belongs to at most $A_2$ discs $D^k$ (we
shall say that $A_2$ is the multiplicity of this covering). Moreover,
\begin{equation}
\label{caracdec}
\eta(D^k)>A_3PT_{1/f}(1)r_k.
\end{equation}
Here $A_1$, $A_2$ and $A_3$ are absolute positive constants and $\eta$
is the non-negative measure on the closed unit disc $\overline\Dbb$
defined by \eqref{caracE}.
\end{theorem} 

\begin{proof*}{\it Proof of Lemma \ref{slowdec}.} Let
\begin{multline*}
\mathcal M_0=\{m: N_m\dm>1/2\},\\
\mathcal M_j=\{m: 2^{-j-1}<N_m\dm\le 2^{-j}\},\ j=1,2,\dots,\\
G_j=\bigcup_{m\in\mathcal M_j}Y_m,\quad P_j=c(\d)2^j,
\end{multline*}
where $c(\d)$ is a positive constant which will be chosen later,
depending only on $\d$, and $\d$ is the constant from \eqref{eqsepar}. Let a
function $f\in H^\infty$, $f\not\equiv0$, be given and let
$D_j^k(w_{j,k},r_{j,k})$ be a system of disks from Theorem~\ref{GovGrish} with
$P=P_j$. We introduce the following sets of disks:
\begin{eqnarray*}
\mathcal D_{j,1}&=&\{D_j^k:D_j^k\cap G_j\ne\varnothing,\quad
r_{j,k}<\tfrac18(1-|w_{j,k}|)\},\quad \mathcal D_1=\bigcup_j\mathcal D_{j,1};\\
\mathcal D_{j,2}&=&\{D_j^k:D_j^k\cap G_j\ne\varnothing,\quad
r_{j,k}\ge\tfrac18(1-|w_{j,k}|)\},\quad \mathcal D_2=\bigcup_j\mathcal D_{j,2}.
\end{eqnarray*}
We shall show that
$$
\sum_{z_l\in\mathcal D_1}(1-|z_l|)<\infty.
$$
By \eqref{eqsepar}, each disk $D_j^k\in\mathcal D_1$ contains at most $N=N(\d)$
points $z_l$. Moreover, disks from $\mathcal D_1$ do not intersect
$\partial\Dbb$, and by \eqref{caracE}
$$
\eta(D_j^k)=\sum_{a_i\in D_j^k}\log\frac1{|a_i|},\quad D_j^k\in\mathcal D_1,
$$
where $a_i$ are zeros of $f$. By \eqref{caracdec}, $\eta(D_j^k)>0$. Hence,
every disc $D_j^k\in\mathcal D_1$
contains at least one zero of $f$. Each disk in $\mathcal D_1$ intersects
at most two sets $Y_m$. For each $j$ the multiplicity of the family
$\{D_j^k\}$ is the absolute constant $A_2$. Hence, the multiplicity of the
covering $\mathcal D_1$ is at most $3A_2$. Thus,
\begin{multline*}
\sum_{z_l\in\mathcal D_1}(1-|z_l|)\le
\sum_{D_j^k\in\mathcal D_1}\sum_{z_l\in D_j^k}(1-|z_l|)\le\\
\le c\sum_{D_j^k\in\mathcal D_1}\sum_{a_i\in D_j^k}(1-|a_i|)
\le 3A_2c\sum_i(1-|a_i|)<\infty,
\end{multline*}
that what we need.

Let $\mathcal J_1$ be the set of indices $m$ such that
$\#\{z_l:z_l\in Y_m\cap\mathcal D_1\}\ge N_m/6$. Then
$$
\sum_{m\in\mathcal J_1}N_m2^{-m}\le c\sum_{z_l\in\mathcal D_1}(1-|z_l|)<\infty.
$$
We set $\mathcal K_1=\mathcal K\setminus\mathcal J_1$. Then at least $\frac56 N_m$
points $z_l$ do not belong to $\mathcal D_1$ for every $m\in\mathcal K_1$.

Fix $j\ge0$. We split the points $z_l \in \mathcal D_{j,2}\cap G_j$,
into two subsets $Z_{j,1}$ and $Z_{j,2}$ in the following way. Let
$z_l\in Y_m$. Then $z_l\in Z_{j,1}$, if every disk $D_j^k\in\mathcal D_{j,2}$
containing $z_l$, does not contain other points $z_{l'}\in Y_m$;
$z_l\in Z_{j,2}$, if there exists a disk $D_j^k\in\mathcal D_{j,2}$
containing at least one other point $z_{l'}\in Y_m$ besides $z_l$. Let
\begin{multline*}
N_m^{j,i}=\#\{z_l:z_l\in Y_m\cap Z_{j,i}\},\quad i=1,2;\\
\mathcal J_{j,2}=\{m\in\mathcal M_j:N_m^{j,1}\ge N_m/6\},\quad
\mathcal J_2=\bigcup_j\mathcal J_{j,2}.
\end{multline*}
We shall prove that
$$
\sum_{m\in\mathcal J_2}N_m2^{-m}<\infty.
$$
Assume that $D_j^k\in\mathcal D_{j,2}$ and $D_j^k\cap Y_m\ne\varnothing$.
Then
\begin{equation}
\label{eq65}
r_{j,k}>\tfrac19 2^{-m-1}.
\end{equation}
Indeed, if $r_{j,k}\le\tfrac19 2^{-m-1}$, then by the inequality
$|w_{j,k}|<1-2^{-m-1}+r_{j,k}$ we have
$$
r_{j,k}\ge\tfrac18(1-|w_{j,k}|)>\tfrac18(2^{-m-1}-r_{j,k})
\ge\tfrac18\cdot8r_{j,k}=r_{j,k},
$$
and we come to a contradiction. From \eqref{eq65} we deduce that
$$
\sum_{m:D_j^k\cap Y_m\ne\varnothing}2^{-m}
<2\cdot2\cdot9r_{j,k}=36r_{j,k}.
$$
Moreover, the number of disks $D_j^k\in\mathcal D_{j,2}$ intersecting
$Y_m$ is at least $N_m^{j,1}$. Hence,
\begin{multline*}
\sum_{m\in\mathcal J_{j,2}}N_m2^{-m}\le
6\sum_{m\in\mathcal J_{j,2}}N_m^{j,1}2^{-m}\\
\le 6 \sum_{k:D_j^k\in\mathcal D_{j,2}}
\sum_{m:D_j^k\cap Y_m\ne\varnothing}2^{-m}
<6\cdot36\sum_k r_{j,k}<\frac c{P_j}=c_12^{-j}
\end{multline*}
(in the last inequality we used \eqref{sumradii}). We have
$$
\sum_{m\in\mathcal J_2}N_m2^{-m}
=\sum_j\sum_{m\in\mathcal J_{j,2}}N_m2^{-m}
<c\sum_j2^{-j}<\infty.
$$
Now we shall show that
\begin{equation}
\label{eq66}
N_m^{j,2}\le\tfrac16N_m \text{ for all } m\in\mathcal K\cap\mathcal M_j.
\end{equation}
For fixed $m$ we consider disks $D_j^k$ containing more than one point
$z_l\in Y_m\cap Z_{j,2}$. Then $2r_{j,k} > d_{m,l}$.
The separation condition \eqref{eqsepar} implies that $d_{m,l} \ge c(\delta) 2^{-m}$,
where $\d$ is the constant from \eqref{eqsepar};
and since the disks $D(z_l, d_{m,l}/2)$ are disjoint, denoting by
$\lambda_2$ the area measure,
\begin{eqnarray*}
C(\delta) 2^{-m} r_{j,k} 
&\ge& \lambda_2 (D(w_{j,k},2r_{j,k})\cap Y_m)
\\
&\ge&
\sum_{l:z_l\in D_j^k\cap Y_m\cap Z_{j,2}}
\lambda_2 \left(D(z_l, d_{m,l}/2)\cap Y_m \right)
\\
&\ge& c'(\delta) \,  2^{-m} \sum_{l:z_l\in D_j^k\cap Y_m\cap Z_{j,2}}
d_{m,l}.
\end{eqnarray*}
Thus
$$
\sum_{l:z_l\in D_j^k\cap Y_m\cap Z_{j,2}}d_{m,l}<c_2(\d)r_{j,k}.
$$
If \eqref{eq66} does not hold, then
\begin{multline*}
\frac16N_m\dm<N_m^{j,2}\dm
\le\sum_{l:z_l\in Y_m\cap Z_{j,2}}d_{m,l}\\
\le\sum_k\sum_{l:z_l\in D_j^k\cap Y_m\cap Z_{j,2}}d_{m,l}
\le\sum_k c_2(\d)r_{j,k}<\frac{c_2(\d)}{P_j}
=\frac{c_2(\d)}{c(\d)}2^{-j}<\frac16N_m\dm,
\end{multline*}
if $c(\d)$ is sufficiently big. That contradiction proves
the validity of \eqref{eq66}.
Thus, if $\mathcal J=\mathcal J_1\cup\mathcal J_2$, then for every
$m\in\mathcal K\setminus\mathcal J$ the number of points $z_l$
lying outside of $\bigcup_kD_j^k$ is at least
$$
\tfrac56N_m-N_m^{j,1}-N_m^{j,2}>\tfrac56N_m-\tfrac16N_m
-\tfrac16N_m=\tfrac12N_m.
$$
The inequality \eqref{eqslowdec} follows from \eqref{logfP}. 
\end{proof*}

\vskip1cm
Vladimir Ya. Eiderman, Moscow State Civil Engineering University, 
129337 Moscow, Yaroslavskoe shosse, 26, Russia. 

eiderman@orc.ru

\vskip.3cm

Pascal J. Thomas, Laboratoire de Math\'ematiques Emile Picard, UMR CNRS 5580,
Universit\'e Paul Sabatier, 118 route de Narbonne, 31062 Toulouse CEDEX, France.

pthomas@cict.fr

\end{document}